\documentclass[reqno,11pt]{amsart}  
\usepackage[fontsize=13]{scrextend}  

\usepackage[hmargin= 25 mm, vmargin= 15 mm]{geometry}

\usepackage{amsmath,amssymb,amsfonts,amsthm}
\usepackage{graphicx}

\usepackage{fmtcount}\usepackage{enumitem}

\newtheorem{theorem}{Theorem}[section]
\newtheorem{lemma}[theorem]{Lemma}

\newtheorem{corollary}[theorem]{Corollary}

\theoremstyle{definition}

\newtheorem{remark}[theorem]{Remark}

\def \ds{\displaystyle}

\newenvironment{alist} {\begin{enumerate}[label=$({\alph*})$,leftmargin=*,labelindent=3mm,labelsep=4mm]} {\end{enumerate}}
\def\bal{\begin{alist}} \def\eal{\end{alist}}

\newenvironment{ilist} {\begin{enumerate}[label=$({\roman*})$,leftmargin=*,labelindent=4mm,labelsep=4mm ]} {\end{enumerate}}
\def\bil{\begin{ilist}} \def\eil{\end{ilist}}

\newenvironment{blist} {\begin{enumerate}[label=\textbullet,leftmargin=*,labelindent=4mm,labelsep=3mm]} {\end{enumerate}}
\def\bbl{\begin{blist}} \def\ebl{\end{blist}}

\newenvironment{blistsmall} {\begin{enumerate}[label=\textbullet,leftmargin=*,labelindent=4mm,labelsep=3mm,noitemsep]} {\end{enumerate}}
\def\bbls{\begin{blistsmall}} \def\ebls{\end{blistsmall}}

\newcommand\ga{\gamma}

\newcommand\de{\delta}
\newcommand\De{\Delta}
\newcommand\ep{\epsilon}

\newcommand\La{\Lambda}
\newcommand\la{\lambda}
\newcommand\om{\omega}
\newcommand\Om{\Omega}
\newcommand\si{\sigma}
\newcommand\Si{\Sigma}

\newcommand \ela{e_\la}
\newcommand \sila{\si_\la}
\newcommand \vla{v_{\la g,v_0}}
\newcommand \Tla{T_{\la g,v_0}}

\newcommand \limn{ \lim_{n \to \infty} }
\newcommand \limt{ \lim_{t \to \infty} }

\newcommand\lamin{\la_{\rm min}  }
\newcommand\psimin{\psi_{\rm min}}
\newcommand\lamax{\la_{\rm max}  }
\newcommand\psimax{\psi_{\rm max}}

\newcommand\bOm{\overline \Om}
\newcommand\Wp{W_0^{1,p}(\Om)}
\newcommand\Wpd{W^{-1,p'}(\Om)}
\newcommand\Wpp{W_{0,+}^{1,p}(\Om)}

\newcommand\Cz{C^0(\bOm)}
\newcommand\Czp{C_+^0(\bOm)}

\newcommand\hf{\hat f}
\newcommand\hv{\hat v}

\newcommand\bg{\overline g}

\newcommand\bv{\overline v}

\newcommand\emb{\hookrightarrow}
\newcommand{\tow}{\rightharpoonup}

\let\le=\leqslant
\let\ge= \geqslant

\newcommand\tom{\widetilde \om}

\newcommand\E{{\mathcal E}}
\newcommand\X{\times}

\def\R{{\mathbb R}}
\newcommand{\pa}{\partial}

\newcommand{\toup}{\scalebox{0.7}{$\,\nearrow\,$}}    
\newcommand{\toups}{\scalebox{0.5}{$\,\nearrow\,$}}  
\newcommand{\todowns}{\scalebox{0.5}{$\,\searrow\,$}}  

\renewcommand{\proof}{\medskip\noindent{\em Proof.} \ }
\def \eop {\qed\medskip}

\numberwithin{equation}{section}

\newcommand{\eqn}[1]{\begin{alignat}{10} #1 \end{alignat}}
\newcommand{\eqns}[1]{\begin{alignat*}{10} #1 \end{alignat*}}
\newcommand{\eqnm}[1]{\begin{equation}\begin{alignedat}{10} #1 \end{alignedat}\end{equation}}

\begin{document}

\title[Global asymptotic stability of positive solutions]{Global asymptotic stability of bifurcating, positive equilibria of $p$-Laplacian boundary value problems with $p$-concave nonlinearities}
\author{Bryan P. Rynne}
\address{Department of Mathematics and the Maxwell Institute for Mathematical Sciences, Heriot-Watt University, Edinburgh EH14 4AS, Scotland.}
\email{B.P.Rynne@hw.ac.uk}

\begin{abstract}
We consider the parabolic, initial value problem
\eqns{v_t&=\Delta_p(v)+\la g(x,v)\phi_p(v),&\quad&\text{in $\Om\X(0,\infty),$}\\
v&=0,&&\text{in $\pa\Om\X(0,\infty),$}\tag{IVP}
\\v&=v_0\ge0,&&\text{in $\Om\X\{0\},$}}
where $\Om$ is a bounded domain in  $\R^N$, for some integer $N\ge1$, with smooth boundary $\pa\Om$,
$\phi_p(s):=|s|^{p-1}\mathop{\rm sgn}s$, $s\in\R$,
$\De_p$ denotes the $p$-Laplacian, with $p>\max\{2,N\}$, $v_0\in\Cz$, and $\la>0$.
The function $g:\bOm\X[0,\infty)\to(0,\infty)$ is $C^0$ and, for each $x\in\bOm$, the function $g(x,\cdot):[0,\infty)\to(0,\infty)$ is Lipschitz continuous and strictly decreasing.

Clearly, (IVP) has the trivial solution  $v\equiv0$, for all $\la>0$.
In addition, there exists $0<\lamin(g)<\lamax(g)$
($\lamax(g)$ may be $\infty$)
such that:
\bbls\item
if $\la\not\in(\lamin(g),\lamax(g))$ then (IVP) has no non-trivial, positive equilibrium;
\item
if $\la\in(\lamin(g),\lamax(g))$ then (IVP) has a unique, non-trivial, positive equilibrium $\ela\in\Wp$.\ebls
We prove the following results on the positive solutions of (IVP):
\bbls\item
if $0<\la<\lamin(g)$ then the trivial solution is globally asymptotically stable;
\item
if $\lamin(g)<\la<\lamax(g)$ then $\ela$ is globally asymptotically stable;
\item
if $\lamax(g)<\la$ then any non-trivial solution blows up in finite time.\ebls
\end{abstract}

\maketitle

\section{Introduction}  \label{intro.sec}

We consider the parabolic, initial-boundary value problem
\eqnm{  \label{orig_t_dept_ivp.eq}
v_t &= \Delta_p(v) + \la g(x,v) \phi_p(v) ,
&\quad& \text{in $\Om \X (0,\infty),$}
\\
v &= 0 ,   && \text{in $\pa \Om \X (0,\infty),$}
\\
v &= v_0 \ge 0 ,  && \text{in $\Om \X \{0\},$}
}
where $\Om$ is a bounded domain in  $\R^N$, for some integer $N \ge 1 $, with smooth boundary $\pa \Om$,
$\phi_p(s)  :=  |s|^{p-1} \mathop{\rm sgn} s$, $s \in \R$,
and $\De_p$ denotes the $p$-Laplacian,
with $p > \max\{ 2,N \} $,
$ v_0 \in \Cz $,
and  $\la > 0$.

We suppose that
$g : \bOm \X [0,\infty) \to (0,\infty) $ is $C^0$
and, for each $x \in \bOm,$
\eqn{  \label{g_dec_strict.eq}
g(x,\cdot) & : [0,\infty) \to (0,\infty) \ \text{is strictly decreasing,}
\\
\label{gainfty_def.eq}
0 \le g_\infty(x) & := \lim_{\xi\to\infty} g(x,\xi) < g_0(x) :=  g(x,0) .
}
We also suppose that $g$ is Lipschitz with respect to $\xi$, in the following sense:
for any $K > 0$ there exists $L_K$ such that
\eqn{  \label{g_lip.eq}
| g(x, \xi_1) - g(x,\xi_2) | \le L_K |  \xi_1 - \xi_2 | ,
  \quad  x \in \bOm , \ 0 \le \xi_1 ,\, \xi_2 \le K .
}

We are interested in positive solutions of \eqref{orig_t_dept_ivp.eq}, so we introduce the following notation:
$ \Czp $ (respectively $\Wpp$) denotes the set of $ \om \in \Cz $ (respectively $ \om \in \Wp$) with $\om \ge 0$ on $\Om$.

It is known that for any $v_0 \in \Czp $ the problem \eqref{orig_t_dept_ivp.eq} has a unique, positive solution
$t \to \vla(t) \in \Wpp $, on some maximal interval $(0,T)$, where we may have  $T < \infty$ or $T=\infty$
(what we mean by a solution will be made precise in Theorem~\ref{existence_time_dept_soln.thm} below).
We are interested in the asymptotic behaviour of these solutions.
This asymptotic behaviour is determined by the  structure of the set of positive equilibria of \eqref{orig_t_dept_ivp.eq}, so we first describe this.

For a given $\la > 0$, a positive {\em equilibrium} is a time-independent solution $u \in \Wpp $ of \eqref{orig_t_dept_ivp.eq},
that is, $u$ satisfies
$\De_p(u) + \la g(u) \phi_p(u)$ = 0
(this will be made precise in Section~\ref{exist-equbsolns.sec} below).
For any $\la > 0$ the function $v \equiv 0$ (or $(\la,v) = (\la,0)$) is a ({\em trivial}) equilibrium.
In addition, the complete structure of the set of non-trivial, positive equilibria of \eqref{orig_t_dept_ivp.eq} is as follows
(see  Theorem~\ref{un_sol_ex.thm} below).
There exists $0 < \lamin(g) < \lamax(g) $
(we may have $ \lamax(g) = \infty $)
such that:
\bbl
\item
if  $ \la \not\in (\lamin(g),\lamax(g)) $ then \eqref{orig_t_dept_ivp.eq} has no non-trivial equilibrium in $\Wpp$;
\item
if $ \la \in (\lamin(g),\lamax(g)) $ then \eqref{orig_t_dept_ivp.eq} has a unique, non-trivial equilibrium $\ela \in \Wpp$.
\ebl

We will prove the following results on the asymptotic behaviour of the positive solutions of \eqref{orig_t_dept_ivp.eq}.
For any $ 0 \ne v_0 \in \Czp $:
\bbls
\item
$ 0 < \la < \lamin(g) \implies \ds\lim_{t \to \infty} \| \vla(t) \|_{0,p} = 0$
\item
$ \lamin(g) < \la < \lamax(g) \implies \ds\lim_{t \to \infty}  \| \vla(t) - \ela \|_{0,p} = 0 $
\item
$ \lamax(g) < \la \implies $ there exists $T < \infty$ such that
$\ds\lim_{t \toups T} |\vla(\cdot)|_0 = \infty$
\ebls

Regarding \eqref{orig_t_dept_ivp.eq} as a bifurcation problem, these results can be interpreted as saying that:
\bbls
\item
when $\la < \lamin(g)$, the trivial equilibrium is globally stable;
\item
as $\la$ increases through $\lamin(g)$, the solution $(\la,0)$ loses stability, and a continuum, $\E^+$, of globally stable, positive equilibrium solutions bifurcates from the point $(\lamin(g),0)$
\\
(in a sense, there is a supercritical, transcritical bifurcation at $\lamin(g)$, with exchange of stability between the equilibria);
\item
as $\la$ increases through $\lamax(g)$, the continuum $\E^+$ `meets infinity' and then disappears, after which all non-trivial, positive solutions blow up in finite time.
\ebls

These results are consistent with a bifurcation analysis of the corresponding semilinear ($p=2$) problem, using the `principle of linearised stability' to obtain local stability.
Such problems have been extensively investigated, see \cite{KS1} and the references therein for a summary of the main results.
However, we do not use bifurcation theory to obtain our results, which usually yields local stability results.
Instead, we use a mixture of comparison and compactness arguments to obtain the above results.

For the quasilinear problem involving the $p$-Laplacian with $p>2$ considered here, these results are consistent with the results on `linearised stability' in the `$p$-concave' case in \cite{KS2}
(condition \eqref{g_dec_strict.eq} is termed `$p$-concavity' in \cite{KS2}; this terminology has been used in other publication for very similar, but slightly different, conditions).
However, the term `linearised stability' in \cite{KS2} refers to the sign of the principal eigenvalue of the linearisation of the problem at an equilibrium solution $\ela$, not to the dynamic (time-dependent) stability that we consider.
In the quasilinear case it is not clear that `linearised stability', in this sense, implies stability in the usual dynamic sense.
Even if such a result could be proved, it would give local rather than global stability.

The convergence results that we obtain say nothing about the rate of convergence.
In particular, we do not obtain the exponential convergence that would be obtained from any sort of `linearised stability' analysis, if such were possible.
Convergence rates for quasilinear problems are discussed in \cite{CF}, together with a broad survey of the literature relating to this.
It is also noted in \cite{CF} that the known results are limited, and difficult to apply.
In particular, the results discussed in \cite{CF} say nothing about the problem considered here.

\section{Preliminaries}  \label{Preliminaries.sec}

\subsection{Notation} \label{Notation.sec}
We let $C^0(\bOm)$ denote the standard space of real valued, continuous functions defined on $\bOm$, with the standard sup-norm on $|\cdot|_0$
(throughout, all function spaces will be real);
$L^q(\Om)$, $q > 1$, denotes the standard space of functions on $\Om$ whose $q$th power is integrable, with norm $\|\cdot\|_q$;
$\Wp$ denotes the standard, first order Sobolev space of functions on $\bOm$ which are zero on $\pa \Om$, with norm  $\|\cdot\|_{1,p}$.
By our assumption that $p > N$, the space $\Wp$ is compactly embedded into $C^0(\bOm)$.
We also define the set of {\em positive} functions in $\Wp$ to be
$\Wpp := \{ \om \in \Wp : \text{$\om \ge 0$ on $\Om$}\}$.
The dual space of $\Wp$ is denoted by $ \Wpd $,
where $p' := p/(p-1)$ is the conjugate exponent of $p$.

If $h : \bOm \X [0,\infty) \to \R$ is continuous then,
for any $\om \in \Czp $, we define $h(\om) \in \Czp$ by
$$
h(\om)(x) := h(x,\om(x)), \quad x \in \bOm .
$$
Clearly, the `Nemitskii' mapping
$\om \to h(\om) : \Czp \to \Czp$ is continuous.
In particular, we repeatedly use the Nemitskii mapping
$\phi_p : \om \to  \phi_p(\om) : \Czp \to \Czp $.

\subsection{The $p$-Laplacian} \label{p-Laplacian.sec}
Formally, the $p$-Laplacian is defined by
$$
\De_p \om := \nabla \cdot (| \nabla \om |^{p-2} \nabla \om) ,
$$
for suitable $\om$,
where $ |{\boldsymbol v}| := (v_1^2 + \dots + v_N^2 )^{1/2} $ for $ {\boldsymbol v} \in \R^N$.
More precisely, for any $\om \in \Wp $, we define
$ \Delta_p(\om) \in  \Wpd $ by
\eqn{   \label{plap_f_weak.eq}
\int_\Om  \Delta_p(\om) \, \varphi
:=
- \int_\Om | \nabla \om |^{p-2} \nabla \om \cdot \nabla \varphi ,
\quad \forall \varphi \in \Wp .
}
A precise definition of what is meant by a solution of \eqref{orig_t_dept_ivp.eq}
will be given in Section~\ref{Time-dependent_solutions.sec} below.

\subsection{Principal eigenvalues of the $p$-Laplacian} \label{evals.sec}

We briefly consider the weighted, nonlinear eigenvalue problem
\eqnm{ \label{princ_eval.eq}
 -\Delta_p(\psi) &= \mu \rho \phi_p(\psi) , \quad  \psi \in \Wp ,
}
where $\mu \in \R$ and the weight function $\rho \in L^1(\Om)$.
We say that $\mu$ is an  {\em eigenvalue} of \eqref{princ_eval.eq}, with {\em eigenfunction} $\psi \in \Wp \setminus \{ 0 \}$,
if the following weak formulation of \eqref{princ_eval.eq} holds
\eqn{   \label{princ_eval_weak.eq}
\int_\Om | \nabla \psi |^{p-2} \nabla \psi \cdot \nabla \varphi
=
\mu  \int_\Om \rho \phi_p(\psi)  \varphi ,
\quad \forall \varphi \in \Wp .
}
A {\em principal eigenvalue} of  \eqref{princ_eval.eq} is an eigenvalue
$\mu_0$  which has a positive eigenfunction $\psi_0 \in \Wp$
(which we will normalise by, say, $|\psi_0|_0 = 1$).
The following result is well known --- see, for example,
\cite[Sections 3-4]{CUE}.

\begin{lemma}  \label{princ_eval.lem}
Suppose that the weight function $\rho$ satisfies: $\rho \ge 0$ on $\Om$,
with $\rho > 0$ on a set of positive Lebesgue measure.
Then the eigenvalue problem \eqref{princ_eval.eq} has a unique principal eigenvalue $\mu_0(\rho)$.
This eigenvalue has the properties, $\mu_0(\rho) > 0$, $\psi_0(\rho) > 0$ on $\Om$,
and
\eqn{   \label{princ_eval_poincare.eq}
\int_\Om |\nabla \om |^p  \ge \mu_0(\rho) \int_\Om \rho |\om|^p  ,
\quad  \forall\, \om \in \Wp .
}
In addition, if $\rho_1,\,\rho_2$ are two such weight functions,
then
\[
\text{$\rho_1 \le \rho_2$ on $\Om$ and
$\rho_1 < \rho_2$ on a set of positive Lebesgue measure}
\implies  \mu_0(\rho_1) > \mu_0(\rho_2).
\]
\end{lemma}

Hence, by \eqref{gainfty_def.eq} and  Lemma~\ref{princ_eval.lem}, we may define
\eqns{
0 < \lamin(g) &:= \mu_0(g_0) < \lamax(g) \ &:=
\begin{cases}
\mu_0(g_\infty) < \infty , & \text{if $g_\infty \ne 0$ \ (in
$L^\infty(\Om)$),}
\\
\infty ,                    & \text{if $g_\infty = 0$ \ (in
$L^\infty(\Om)$).}
\end{cases}
}
and we denote the corresponding normalised principal eigenfunctions by $\psimin(g),\, \psimax(g)$.

\subsection{An energy functional}  \label{energy_fnl.sec}

We now define an `energy' functional for \eqref{orig_t_dept_ivp.eq} on $\Wpp$.
Let
\eqns{
F(x,\xi) &:= \int_0^\xi g(x,s) s^{p-1}\,ds,
&\quad&  (x,\xi) \in \Om \times [0,\infty),
\\
E_{\la g}(\om) &:= \frac{1}{p} \int_\Om |\nabla \om|^p - \la \int_\Om F(\om),
 && \om \in \Wpp .
}
By the continuity of the embedding $\Wpp \emb C^0(\bOm)$,
the {\em energy} functional $E_{\la g} : \Wpp \to \R$ is continuous.

\begin{lemma}  \label{Ebdd.lem}
If $\la < \mu_0(g_\infty) $ then there exists an increasing function $M_\la : \R \to (0,\infty)$ such that,
$$
|\om|_0  + \|\om\|_{1,p} < M_\la(E_{\la g}(\om)) ,
\quad \om \in \Wpp .
$$
\end{lemma}

\proof
Suppose the contrary, so there exists $R \in \R$ and
$0 \ne \om_n \in \Wpp$, $n = 1,2,\dots,$
such that
$E_{\la g}(\om_n) \le R$ and $\lim_{n \to \infty} \|\om_n\|_{1,p} = \infty$
(since $p > N$, $|\om|_0 \le C_0 \|\om\|_{1,p}$, for some constant $C_0$).
Let $\tom_n := \om_n / \|\om_n\|_{1,p}$, $n = 1,2,\dots$.
By the compactness of the embedding $\Wp \emb C^0(\bOm)$,
we may assume that $\tom_n \to \tom_\infty $ in $\Czp$, for some
$\tom_\infty \in \Czp$,
and it suffices to show that this leads to a contradiction.

By definition,
\eqn{  \label{Ewn.eq}
E_{\la g}(\om_n) = \frac{1}{p} \| \om_n \|_{1,p}^p
\left\{ \int_\Om |\nabla \tom_n|^p - \la p \int_\Om \frac{F(\om_n)}{\|\om_n\|_{1,p}^p} \right\} ,
\quad n \ge 1 .
}
We now show that, as $n \to \infty$,
\begin{equation}  \label{F_lim.eq}
p \int_\Om  \frac{F(\om_n)}{\|\om_n\|_{1,p}^p}
\to
\int_\Om g_\infty \tom_\infty^p .
\end{equation}
By \eqref{g_dec_strict.eq} and \eqref{gainfty_def.eq}
there exists $C > 0$ such that, for any  $n \ge 1$,
\eqn{  \label{unifbndF.eq}
p \frac{|F(\om_n)|_0}{\|\om_n\|_{1,p}^p}
\le
\frac{|g_0|_0 |\om_n|_0^p}{\|\om_n\|_{1,p}^p}
\le  C ,
}
and similarly, using \eqref{gainfty_def.eq}, for any $x \in \Om$ and  $\ep > 0$, there exists $C(x,\epsilon) > 0$ such that, for any  $n \ge 1$,
$$
p \frac{F(\om_n)(x)}{\|\om_n\|_{1,p}^p}
\le
\frac{C(x,\epsilon) + (g_\infty(x) + \epsilon)
\om_n(x)^p}{\|\om_n\|_{1,p}^p}
\to (g_\infty(x) + \epsilon) \tom_\infty(x)^p  .
$$
Combining this with a similar lower bound shows that
\eqn{  \label{Frescale_lim.eq}
p \frac{F(\om_n)(x)}{\|\om_n\|_{1,p}^p}
\to g_\infty(x) \tom_\infty(x)^p  ,
\quad  x \in \Om ,
}
so \eqref{F_lim.eq} follows from \eqref{unifbndF.eq}, \eqref{Frescale_lim.eq} and the dominated convergence theorem.

Now suppose that $\int_\Om g_\infty \tom_\infty^p > 0$.
Then, by Lemma~\ref{princ_eval.lem}, for $n \ge 1$,
\begin{equation}  \label{est_pos.eq}
\int_\Om |\nabla \tom_n |^p
\ge \mu_0(g_\infty)  \int_\Om g_\infty \tom_n^p
\to  \mu_0(g_\infty)  \int_\Om g_\infty \tom_\infty^p
> 0,
\end{equation}
and combining \eqref{Ewn.eq}, \eqref{F_lim.eq} and \eqref{est_pos.eq} shows that $E_{\la g}(\om_n) \to \infty$
(since $\la < \mu_0(g_\infty)  $).
However, this contradicts the initial assumption that
$E_{\la g}(\om_n) \le R$ for all $n \ge 1.$

Next, suppose that $\int_\Om g_\infty \tom_\infty^p = 0$, with $ \|\tom_\infty\|_p > 0. $
Then,  by Lemma~\ref{princ_eval.lem}, for $n \ge 1$,
\begin{equation}  \label{est_zero.eq}
\int_\Om |\nabla \tom_n |^p
\ge \mu_0({\bf 1}) \|\tom_n\|_p^p
\to  \mu_0({\bf 1}) \|\tom_\infty\|_p^p
> 0
\end{equation}
(where $\bf 1$ denotes the weight function that is identically 1 on
$\Om$), and combining \eqref{Ewn.eq}, \eqref{F_lim.eq} and \eqref{est_zero.eq} again yields the contradiction $E_{\la g}(\om_n) \to \infty$.

Finally, suppose that $\|\tom_\infty\|_p = 0$.
Since $\|\tom_n\|_{1,p} = 1$, $n \ge 1$, this implies that
$
\int_\Om |\nabla \tom_n |^p \to 1,
$
so  combining \eqref{Ewn.eq} and \eqref{F_lim.eq} again yields a contradiction, and so  completes the proof of
Lemma~\ref{Ebdd.lem}.
\eop

\subsection{Existence and uniqueness of non-trivial, positive equilibria}  \label{exist-equbsolns.sec}

A {\em positive equilibrium} of \eqref{orig_t_dept_ivp.eq} is a solution of the problem
\eqnm{ \label{equb.eq}
- \De_p(u) &= \la g(u) \phi_p(u) , \quad  u \in \Wpp .
}
More precisely, a {\em solution} of \eqref{equb.eq} is defined to be a function $u \in \Wpp$
which satisfies the following weak formulation of \eqref{equb.eq},
\eqn{   \label{equb_weak.eq}
\int_\Om | \nabla u |^{p-2} \nabla u \cdot \nabla \varphi
=
\la  \int_\Om g(u)  \phi_p(u) \varphi ,
\quad \forall \varphi \in \Wp .
}
Clearly, for any $\la \in \R$, the function $u=0$ is a ({\em trivial}) positive solution of
\eqref{orig_t_dept_ivp.eq} and \eqref{equb.eq}.

We now describe the structure of the set of non-trivial, positive equilibria.
Let
$$
 \La := (\lamin(g),\lamax(g)) .
$$

\begin{theorem} \label{un_sol_ex.thm}
\bal
\item
If $ \la \not\in \La $ then \eqref{equb.eq} has no non-trivial solution $u \in \Wpp$.
\item
If $ \la \in \La $ then \eqref{equb.eq} has a unique, non-trivial solution $\ela \in \Wpp$, and $\ela > 0$ on $ \Om$.
\item
The mapping $\la \to \ela : \La \to \Wpp $ is continuous, and
\eqn{  \label{E_limits.eq}
\lim_{\la \todowns \lamin(g)} \| \ela \|_{1,p} = 0 ,
\quad
\lamax(g) < \infty \implies \lim_{\la \toups \lamax(g)} \| \ela \|_{1,p} = \infty .
}
\eal
\end{theorem}

\proof
Parts~$(a)$ and~$(b)$ are proved in \cite[Theorems~1,~2]{DS}.
We observe that:\\[1ex]
$(i)$ \ the strict positivity of $ \ela $ on $\Om$ is not stated in \cite[Theorem~2]{DS}, but is derived in its proof; it also follows from Lemma~\ref{princ_eval.lem};\\
$(ii)$ \ part~$(a)$ also follows from Lemma~\ref{princ_eval.lem} and the definitions of
$\lamin(g),\,\lamax(g)$.
\medskip

To prove part~$(c)$, suppose firstly that
$\la_n \in \La$, $n=1,2,\dots,$
is such that
\eqn{  \label{rab_sequence.eq}
\lim_{n \to \infty} \la_n = \la_\infty < \infty ,
\quad
\lim_{n \to \infty} \| e_{\la_n} \|_{1,p}  = \infty .
}
Defining
$w_n := e_{\la_n} / \| e_{\la_n} \|_{1,p} $, $n=1,2,\dots,$
it follows from the compactness and continuity properties described on p.~229 of \cite{DM} that we may suppose that $w_n \to w_\infty \in \Wpp$, with $w_\infty \ne 0$ and
\eqnm{  \label{rablim.eq}
 -\Delta_p(w_\infty) &= \la_\infty \bg \phi_p(w_\infty) ,
\\
\bg(x) &= \lim_{n \to \infty} g( e_{\la_n}(x)) , \quad x \in \Om .
}
By \eqref{g_dec_strict.eq} and \eqref{gainfty_def.eq},
$ 0 \le \bg \le g_0 $ on $\bOm$,
so it follows from \eqref{rablim.eq} and the invertibility of the operator $\Delta_p$ (see \cite{DM}) that we must have $ \bg > 0 $ on a set of positive measure.
Hence, by Lemma~\ref{princ_eval.lem} and \eqref{rablim.eq}, $ w_\infty(x) > 0 $ for each  $x \in \Om$, so that $e_{\la_n}(x) \to \infty$, and $\bg(x) = g_\infty(x)$.
Thus, by the definition of $\lamax(g)$ and \eqref{rablim.eq}, $ \la_\infty = \lamax(g) $.
We conclude that the mapping  $\la \to \ela : \La \to \Wpp $ is bounded on any closed, bounded subinterval of
$ \La $, and hence, again using the continuity properties in \cite{DM}, this mapping is continuous on $ \La $.

Next, by similar arguments, it can be shown that if $\la_n \to \lamin(g)$ then $ \| e_{\la_n} \|_{1,p} $ cannot be bounded away from 0,
and if $\la_n \to \lamax(g)$ then $ \| e_{\la_n} \|_{1,p} $ cannot be bounded, which proves \eqref{E_limits.eq},
and so completes the proof of the theorem.
\eop

\remark
$(a)$ \ Theorem~\ref{un_sol_ex.thm} shows that the set of non-trivial, positive equilibria,
which we will denote by $\E^+$,
is a Rabinowitz-type continuum in $ \La \X \Wpp$, which bifurcates from $(\lamin(g),0)$ and `meets infinity' at $\lamax(g)$.
\\[1ex]
$(b)$ \
It is shown in \cite{GS}, \cite{GEN1} that if $\Om$ is a ball, then $\E^+$ is in fact a smooth curve of radially symmetric solutions.

\section{Time-dependent solutions}  \label{Time-dependent_solutions.sec}

In Section~\ref{exist-equbsolns.sec} we discussed equilibrium (time-independent) solutions of equation \eqref{orig_t_dept_ivp.eq}.
In this section we will discuss time-dependent solutions of \eqref{orig_t_dept_ivp.eq}.
We first describe an existence and uniqueness result, and then a comparison result, which will be used to determine the long-time behaviour of the solutions.

\subsection{Existence and uniqueness of positive solutions}  \label{Time-dependent_solutions-existence.sec}

In this section we will discuss the existence, uniqueness and properties of solutions of the time-dependent problem \eqref{orig_t_dept_ivp.eq}.
To state precisely what we mean by a solution of \eqref{orig_t_dept_ivp.eq} we define the spaces
\eqns{
\Si(T) &:= C([0,T),L^2(\Om))
\,\cap\,
C((0,T),\Wp)
\,\cap\,
W_{\rm loc}^{1,2}((0,T),L^2(\Om)) ,
\quad T > 0
}
(we allow $T=\infty$ here, and likewise for other such numbers below).
The space $W^{1,2}((0,T),L^2(\Om)) $ is defined
in \cite[Example~10.2]{RR},
using the notation $H^1((0,T),L^2(\Om))$;
the space $W_{\rm loc}^{1,2}((0,T),L^2(\Om))$ can be defined by a simple adaptation of the definition in \cite{RR}.
We will search for a solution of \eqref{orig_t_dept_ivp.eq} in $\Si(T)$, for some $T > 0$.
Thus, in this setting, a solution $v$ will be regarded as a time-dependent mapping
$ t \to v(t) : (0, T ) \to \Wp $,
with $ \De_p(v(t)) \in \Wpd $ defined  by \eqref{plap_f_weak.eq}, for each $t \in (0, T )$,
and satisfying the initial condition at $ t=0 $ as a limit in $ L^2(\Om) $.
More (or less) regularity at $ t=0 $ can be attained, depending on the regularity of $v_0$
(for example if $v_0 \in \Wpp$ then the solution will belong to $C([0,T),\Wp)$),
but the above setting will suffice here.

In view of this, we will rewrite \eqref{orig_t_dept_ivp.eq} in the form
\eqn{  \label{rw_t_dept_ivp.eq}
\frac{dv}{dt} = \De_p(v) + \la g(v) \phi_p(v), \quad v(0) = v_0 \in \Czp .
}
The following theorem describes the existence and uniqueness of solutions of \eqref{rw_t_dept_ivp.eq}, and various additional properties which will be required below.
This theorem does not require $g$ to be positive, nor to satisfy the conditions \eqref{g_dec_strict.eq}, \eqref{gainfty_def.eq}.
It does, however, assume that $g$ is defined on $\bOm \X \R$ rather than on $\bOm \X [0,\infty)$
(which we have assumed so far, since we are mainly interested in positive solutions).
Once we have established the general existence of solutions,
we will then prove their positivity, and thereafter the values of $g(x,\xi)$, $\xi < 0$, will be irrelevant.
If $g$ is only defined on $\bOm \X [0,\infty)$ ab initio, then we may simply extend it to $\bOm \X \R$ by setting
$g(x,-\xi) = g(x,0)$, $\xi > 0$.

\begin{theorem}  \label{existence_time_dept_soln.thm}
Suppose that $g : \bOm \X \R \to \R $ satisfies the Lipschitz condition \eqref{g_lip.eq}  on $\bOm \X \R$,
and $ \la > 0 $, $  v_0 \in \Cz $.
Then  \eqref{rw_t_dept_ivp.eq}
has a unique solution $\vla \in \Si( \Tla )$,
defined on a maximal interval $[0, \Tla )$,
for some $  \Tla  > 0 $,
having the following properties.
\bal
\item  \label{thmlab-ival_satisfied.item} 
$\vla(0) = v_0.$
\item   \label{thmlab-solve_eqn.item}   
The function
$ \vla : [0, \Tla ) \to L^2(\Om) $ is differentiable
at almost all $t \in [0, \Tla )$, and at such $t$,
$$
\frac{d \thinspace \vla }{dt}\, (t) \, , \ \Delta_p(\vla(t)) \in L^2(\Om) ,
$$
and
$$
\frac{d \thinspace \vla }{dt}\, (t) = \De_p(\vla(t)) + \la g(\vla(t)) \phi_p(\vla(t)) ,
   \quad \text{in $L^2(\Om).$}
$$
\item  \label{thmlab-E_decreasing.item} 
The function
$E_{\la g}(\vla(\cdot)) : (0, \Tla ) \to \R$
is absolutely continuous, decreasing and
\eqn{  \label{E_t.eq}
\frac{d}{dt}\, E_{\la g}(\vla(t)) = -\Big\| \frac{d}{dt}\, \vla(t) \Big\|_2^2 ,
\quad  \text{ a.e. \ $t \in (0, \Tla ) .$ }
}
\item   \label{meaning_maximal.item} 
The interval  $[0, \Tla )$ on which the solution exists is maximal, in the sense that
\eqn{  \label{meaning_maximal.eq}
 \Tla  < \infty  \implies
\limsup_{t \toups  \Tla } | \vla(t) |_0 = \infty .
}
If the set of equilibria of \eqref{rw_t_dept_ivp.eq} is bounded in $\Cz$ then \eqref{meaning_maximal.eq} holds with $\lim$ rather than $\limsup$.
\eal
\end{theorem}

\proof
Let $\theta \in C^\infty(\R,\R)$ be a decreasing function with
$$
\theta(s) =
\begin{cases}
1, & s \le 1 ,
\\
0, & s \ge 2  ,
\end{cases}
$$
and for any integer $n \ge 1$, define
$\hat f_n : \bOm \times \R \to (0,\infty)$ by
$$
\hat f_n(x,\xi) :=
\theta(|\xi|/n)\, g(x,\xi)\, \phi_p(\xi) , \quad (x,\xi) \in \bOm \times \R .
$$
Since $\hat f_n$ is bounded and Lipschitz, the results of \cite{AO} and \cite{CF}  show that the problem
\eqnm{  \label{hv_time_dept.eq}
\hv_t  &= \Delta_p(\hv) + \la \hat f_n(\hv)  ,
\quad \hv(0) &= v_0  ,
}
has a unique solution $\hv_n \in \Si(\infty)$ having the properties
\ref{thmlab-ival_satisfied.item}-\ref{thmlab-E_decreasing.item}
(we discuss this further in  Remark~\ref{strong_soln_2.rem} below).
Clearly, $\hv_n $ is a solution of \eqref{rw_t_dept_ivp.eq} on the time interval $[0,T_n)$, where
$$
T_n := \sup \{ T : |\hv_n(t)|_0 \le n : t \in [0,T) \} , \quad n \ge 1 ,
$$
and letting
$$
 \Tla  := \limn T_n  ,
$$
we see that \eqref{rw_t_dept_ivp.eq} has a unique solution $ \vla \in \Si( \Tla )$,
having the properties
\ref{thmlab-ival_satisfied.item}-\ref{thmlab-E_decreasing.item},
and
\eqn{ \label{meaning_maximal_subsequence.eq}
 \Tla  < \infty  \implies
\lim_{t_n \toups  \Tla } |  \vla (t_n) |_0 = \infty ,
}
for some sequence $(t_n)$ in $(0,\Tla)$.
That is, \eqref{meaning_maximal.eq} holds.

Now suppose
that there exists $M>0$ such that $|e|_0 < M$ for any equilibrium solution $e$ of \eqref{rw_t_dept_ivp.eq},
and that the overall limit in \eqref{meaning_maximal.eq} does not exist.
Then there exists another sequence $(s_n)$ in $(0,\Tla)$ such that $s_n \toup  \Tla$ and the sequence $ (|  \vla (s_n) |_0) $ is bounded.
Combining this with \eqref{meaning_maximal_subsequence.eq} and the continuity of the mapping $ t \to |  \vla (t) |_0 $ on $(0,\Tla)$, we may suppose further that $ | \vla (s_n) |_0 = K $, $n = 1,2,\dots$, for some $K>M$.
It now follows from property~\ref{thmlab-E_decreasing.item} and the definition of $E$ that the sequence $ (\vla (s_n)) $ is bounded in $\Wp$ so,
by the argument in Section~\ref{proof-part_ab1.sec} below
(after taking a subsequence if necessary)
there exists $v_\infty \in \Wp$ such that
$| \vla (s_n) - v_\infty |_0 \to 0$
and $v_\infty$ is an equilibrium of \eqref{rw_t_dept_ivp.eq}
(it is assumed in Section~\ref{proof-part_ab1.sec} that $0 < \la < \lamax(g)$, but this assumption is only used to obtain such a bounded sequence in $\Wp$).
However, this implies that $ |v_\infty |_0 = K > M $, which contradicts the choice of $M$ above, so we conclude that the limit in \eqref{meaning_maximal.eq} must in fact exist, that is, property \ref{meaning_maximal.item} must hold.
This completes the proof of Theorem~\ref{existence_time_dept_soln.thm}.
\eop

\begin{remark}  \label{strong_soln_2.rem}
$(a)$ \
The existence and uniqueness of a solution $\hv_n$ of \eqref{hv_time_dept.eq}, and the fact that $\hv_n$ has properties
\ref{thmlab-ival_satisfied.item}-\ref{thmlab-solve_eqn.item},
as asserted in the proof of Theorem~\ref{existence_time_dept_soln.thm},
follows by combining various standard results on maximal monotone operators.
Specifically,
\cite[Theorem~3.2]{BRE1}, \cite[Theorems~3.4, 3.11]{AO} and  \cite[Remark~3.6(5)]{AO}.
How these results combine to give a solution with the desired properties is discussed in detail in  \cite[Remark~2.2]{CF}.
It should be noted that, with the sign of $f$ used here, the functions $\hf_n$ are not monotone but, by assumption \eqref{g_lip.eq}, they satisfy the Lipschitz condition imposed on the function $f_2$ in assumption $(2.4)$ in \cite{CF}.
Thus, to apply the discussion in \cite{CF} to the problem \eqref{hv_time_dept.eq} above,
we set (in the notation in \cite{CF}) $f_1 = 0$ and $f_2 = \hf_n $.

The fact that $\hv_n$ also has property
\ref{thmlab-E_decreasing.item}
follows from  \cite[Lemma~3.3]{BRE1}, and the argument in the proof of \cite[Lemma~3.1]{CF}.
\\[1ex]
$(b)$ \ Existence and uniqueness of a local (in time) solution of \eqref{rw_t_dept_ivp.eq}, with weaker properties than those stated in Theorem~\ref{existence_time_dept_soln.thm}, is proved in \cite[Theorem~2.1]{JUN},
and global existence and uniqueness of such solutions of \eqref{hv_time_dept.eq}
(under similar Lipschitz conditions)
is proved in \cite[Theorem~3.1]{JUN}.
Hence, the solution $\vla$ given by Theorem~\ref{un_sol_ex.thm} is unique in a considerably broader solution space than $\Si^p$.
\end{remark}

\subsection{Comparison results}

We now consider the auxiliary problem
\eqn{  \label{comparison_gamma.eq}
\frac{d w}{d t}  &= \De_p(w) + \la \ga \phi_p(w) ,  &\quad  w(0) = w_0 \in \Czp ,
}
where $\ga \in L^\infty(\Om)$ is independent of $v$,
and $\ga \ge 0$ on $\Om$.
This is a special case of  \eqref{rw_t_dept_ivp.eq}
(with $g(x,v)$ having the form $\ga(x)$) so, by Theorem~\ref{existence_time_dept_soln.thm}, the problem \eqref{comparison_gamma.eq} has a unique solution $w_{\la\ga,w_0}$ defined on a maximal interval $[0,T_{\la\ga,w_0})$.

\begin{remark}
Theorem~\ref{existence_time_dept_soln.thm} was stated, and proved, for continuous functions $g$ depending on $(x,\xi)$ (and Lipschitz with respect to $\xi$), but the results quoted from \cite{CF}
(see Remark~\ref{strong_soln_2.rem})
in the proof of Theorem~\ref{existence_time_dept_soln.thm} apply equally to the problem  \eqref{comparison_gamma.eq}, containing an $x$-dependent function $\ga \in L^\infty(\Om)$.
\end{remark}

We now describe a `comparison' result for solutions of
\eqref{rw_t_dept_ivp.eq} and \eqref{comparison_gamma.eq}.
For any $T > 0$ and functions $\om_1,\, \om_2 \in \Si(T)$, we write
$ \om_1 \ge \om_2 $ on $[0,T)$ if
$ \om_1(t) \ge \om_2(t) $, on $\bOm$, for each $t \in [0,T)$.
From now on we suppose that $g$ satisfies our basic hypotheses, that is, $g$ is positive and satisfies \eqref{g_dec_strict.eq} and \eqref{gainfty_def.eq}.

\begin{lemma} \label{comp.lem}
If $ g_\infty \ge \ga \ge 0 $ and $v_0 \ge  w_0 \ge 0 $ on $\bOm$, then
$$
 \Tla  \le T_{\la\ga,w_0}
\quad \text{and} \quad
 \vla  \ge w_{\la\ga,w_0}
\quad \text{on $ [0,  \Tla  )$.}
$$
\end{lemma}

\proof
The proof follows, with minor modifications, the proof of \cite[Theorem~2.5]{LX}, using our assumptions on $g$
(in particular, the assumption $g_\infty \ge \ga$ implies that $ g(v) \ge \ga $ for any $v \in \Wpp$).
We omit the details.
However, we note that \cite[Theorem~2.5]{LX} considers equations of the form
$v_t = \De_p(v) + \la \phi_p(v) $, but the proof can be adapted to give the above result;
the argument in \cite{LX} is based on the proof of \cite[Lemma~3.1, Ch.~VI]{DIB}, which considered the equation
$v_t = \De_p(v) $.
\eop

If $\ga = 0$ and $w_0=0$, then clearly $w_{0,0} \equiv 0$, and since $g_\infty \ge 0$, Lemma~\ref{comp.lem} now yields the following positivity result for the solution $\vla$ of \eqref{rw_t_dept_ivp.eq} found in Theorem~\ref{existence_time_dept_soln.thm}.

\begin{corollary}  \label{vla-positive}
If $v_0 \in \Czp $ then $\vla(t)  \in \Wpp  $ for all $t \in (0, \Tla ) .$
\end{corollary}

In the next section we will use the comparison result Lemma~\ref{comp.lem} to describe the behaviour of solutions of \eqref{rw_t_dept_ivp.eq}.
The following criterion for finite time blow-up of solutions of \eqref{comparison_gamma.eq} will be useful.

\begin{lemma} \label{blowup_wga.lem}
If  $ \la > \mu_0(\ga) $ and $ 0 \ne w_0 \in \Czp $, then
$ T_{\la\ga,w_0} < \infty . $
\end{lemma}

\proof
This can be proved by following, almost verbatim, the proof of \cite[Theorem~3.5]{LX}, which deals with the case $\ga \equiv 1$.
\eop

\section{Global stability and instability of the equilibria}  \label{stability.sec}

For any $\la > 0$ the time-dependent problem \eqref{rw_t_dept_ivp.eq}
has the trivial equilibrium solution $ u = 0$,
and also, by Theorem~\ref{un_sol_ex.thm},
for each $\la \in (\lamin(g) ,\lamax(g) )$
there is a unique, non-trivial, positive equilibrium $\ela \in \Wpp$.
We will now consider the global stability, and instability, of these equilibria.

\begin{theorem}  \label{global_asymptotes.thm}
Suppose that $ 0 \ne v_0 \in \Czp .  $
\bal
\item
If $ 0 < \la \le \lamin(g) $  then
$ \Tla  = \infty$ \, and \
$\ds\lim_{t \to \infty} \|\vla(t)\|_{1,p} = 0 . $
\item
If $ \lamin(g) <  \la <  \lamax(g) $  then
$ \Tla  = \infty$ \, and \
$\ds\lim_{t \to \infty} \|\vla(t) - \ela\|_{1,p} = 0 . $
\item
If $ \lamax(g) < \la $  then
$ \Tla  < \infty$,
that is, the solution $ \vla $ blows up in finite time.
\eal
\end{theorem}

\subsection{Proof of Theorem~\ref{global_asymptotes.thm}~$\boldsymbol{(a),\,(b)}$ }
\label{proof-part_ab1.sec}
Suppose that  $ 0 < \la < \lamax(g) = \mu_0(g_\infty) $.
Let $ \bv = \vla(\Tla/2) \in \Wpp$.
Then $E_{\la g}(\bv)$ is defined and, by Lemma~\ref{Ebdd.lem} and Theorem~\ref{existence_time_dept_soln.thm}~$(d)$-$(e)$,
\eqnm{ \label{v_le_MEv.eq}
E_{\la g}(\vla(t)) \le E_{\la g}(\bv)
&\implies
\|\vla(t)\|_{1,p} \le M_\la(E_{\la g}(\bv)) , \quad  \tfrac12 \Tla \le t <  \Tla
\\[1ex]
&\implies
\text{$\Tla  = \infty$ \ and \ $E_{\la g}(\vla(\cdot))$ is bounded on $(0,\infty)$}
\\[1ex]
&\implies
\text{$\ds\lim_{t \to \infty} E_{\la g}(\vla(t))$
exists.}
}

From now on, $(t_n)$ will denote an increasing sequence in $(0,\infty)$ such that $ t_n \to \infty$; we will choose various such sequences below.
By  \eqref{v_le_MEv.eq}, the sequence $(\vla(t_n))$ is bounded in $\Wp$, so we may also suppose (after taking a subsequence if necessary) that
\eqn{ \label{cvgce_vinf_in_wpw.eq}
\vla(t_n) \tow v_\infty \ \text{in $\Wp$} , \quad
|\vla(t_n) - v_\infty|_0 \to 0 ,
}
for some $ v_\infty \in \Wpp $
(where $ \tow$ denotes weak convergence in $\Wp$).
The argument in the proof of \cite[Lemma~3.1]{CF} now shows that $v_\infty$ is an equilibrium solution of \eqref{rw_t_dept_ivp.eq},
that is, $v_\infty$ is a solution of \eqref{equb.eq}.
Hence, by Theorem~\ref{un_sol_ex.thm}, we have the following cases.
\bal
\item
If  $ 0 < \la \le \lamin(g) $  then $v_\infty = 0$.
\item
If  $ \lamin(g) <  \la <  \lamax(g) $ then either  $v_\infty = 0$ or $v_\infty = \ela $.
\eal
In case~$(a)$, a simple contradiction argument (using the preceding results) now shows that we must have $ \limt |\vla(t) - \sila|_0 = 0$.
In case~$(b)$, suppose that there exists sequences $(t_n^0)$, $(t_n^1)$, such that
$$
|\vla(t_n^0)|_0 \to 0 , \quad |\vla(t_n^1) - \ela |_0 \to 0 .
$$
Then, by continuity of the mapping $t \to |\vla(t)|_0$, there exists a sequence $(t_n^2)$ such that
$$
|\vla(t_n^2)|_0 =  | \ela |_0 /2 ,  \quad  n \ge 1 .
$$
But this contradicts the preceding results, so we must have
$ \limt |\vla(t) - v_\infty|_0 = 0$, for either  $v_\infty = 0$ or $v_\infty = \ela $.
The following lemma shows that, in fact, the latter holds -- the proof will use some results from Section~\ref{proof-part_c.sec}, so will be given in Section~\ref{proof-lemma.sec}  below.

\begin{lemma}  \label{vinfty_ne_0.lem}
For $ \lamin(g) < \la < \lamax(g) $ and $ 0 \ne v_0 \in \Czp $, \
$ \ds \limt |\vla(t) -  \ela |_0 = 0$.
\end{lemma}

Now, to simplify the notation, and to combine the two cases $(a)$ and $(b)$, we define $\sila \in \Wpp$ by
$$
\sila :=
\begin{cases}
0 ,  & \text{if $0 < \la \le \lamin(g) $,}
\\
\ela ,  & \text{if $ \lamin(g) < \la < \lamax(g) , $}
\end{cases}
$$
and the preceding results show that
\eqn{  \label{limt_vla.eq}
\lim_{t \to \infty} |\vla(t) - \sila |_0 = 0 , \quad 0 < \la < \lamax(g) .
}
Thus, it only remains to prove the convergence with respect to the $\Wp$ norm.

By integrating \eqref{E_t.eq} with respect to $t$ and using the existence of the limit $\lim_{t \to \infty} E_{\la g}(\vla(t))$
(by \eqref{v_le_MEv.eq}),
we see that the function on the right hand side of \eqref{E_t.eq} lies in $L^2(0,\infty)$,
so we may choose a sequence $(t_n)$ such that
\eqn{ \label{rhs_cvgce_in_L2.eq}
&\| \De_p(\vla(t_n)) + \la g(\vla(t_n)) \phi_p(\vla(t_n)) \|_2 \to 0 .
}
We may also suppose that \eqref{cvgce_vinf_in_wpw.eq} holds, with
$ v_\infty = \sila $.
Hence, by \eqref{plap_f_weak.eq}, \eqref{equb_weak.eq} and \eqref{rhs_cvgce_in_L2.eq},
\eqns{  
&\int_\Om \big(\De_p(\vla(t_n)) + \la g(\vla(t_n)) \phi_p(\vla(t_n)) \big) \vla(t_n) \to 0 ,
&\quad&
\text{(by \eqref{v_le_MEv.eq} and \eqref{rhs_cvgce_in_L2.eq})}
\\[1ex]
&  \implies
\int_\Om |\nabla \vla(t_n)|^p \to \la \int_\Om g(\sila) \sila^p
=  \int_\Om |\nabla \sila|^p ,
&\quad&
\text{(by \eqref{plap_f_weak.eq} and \eqref{equb_weak.eq})}
}
and combining this with \eqref{cvgce_vinf_in_wpw.eq} yields
\eqn{  \label{wk_and_nm_cvgce.eq}
\vla(t_n) \tow \sila, \quad \ \text{in $\Wp$} , \quad
\|\vla(t_n)\|_{1,p} \to \| \sila \|_{1,p}  .
}
Hence, by the uniform convexity of $\Wp$ and \cite[Proposition 3.32]{BRE2},
$
\|\vla(t_n) - \sila \|_{1,p}  \to 0 ,
$
which implies that
$
E_{\la g}(\vla(t_n)) \to E_{\la g}(\sila) ,
$
and so
\eqn{  \label{Elim.eq}
\lim_{t \to \infty} E_{\la g}(\vla(t)) = E_{\la g}(\sila)
}
(since this limit exists, by \eqref{v_le_MEv.eq}).

Now suppose that there exists a sequence $(t_n)$ and $\ep > 0$ such that
$ \|\vla(t_n) - \sila \|_{1,p}  > \ep  $, and also that \eqref{cvgce_vinf_in_wpw.eq} holds.
Combining this with \eqref{Elim.eq}
(and the form of $ E_{\la g} $)
shows that \eqref{wk_and_nm_cvgce.eq} again holds, and so (by uniform convexity)
$
\|\vla(t_n) - \sila \|_{1,p}  \to 0 ,
$
which contradicts this choice of sequence $(t_n)$.
Hence, we must have
$
\|\vla(t) - \sila \|_{1,p}  \to 0 ,
$
which completes the proof of parts~$(a)$ and~$(b)$ of Theorem~\ref{global_asymptotes.thm}.

\subsection{Proof of Theorem~\ref{global_asymptotes.thm}~$\boldsymbol{(c)}$ }
\label{proof-part_c.sec}
By hypothesis, $\la > \lamax(g) = \mu_0(g_\infty)$,
so by Lemmas~\ref{comp.lem} and~\ref{blowup_wga.lem}
(with $\ga = g_\infty$ and $w_0 = v_0$),
$
\Tla \le T_{\la g_\infty, v_0} < \infty ,
$
which proves part~$(c)$ of Theorem~\ref{global_asymptotes.thm}.
\eop

\subsection{Proof of Lemma~\ref{vinfty_ne_0.lem} }
\label{proof-lemma.sec}

By the arguments preceding Lemma~\ref{vinfty_ne_0.lem} in Section~\ref{proof-part_ab1.sec}, it suffices to show that if we suppose that
\eqn{  \label{vla_to_z.eq}
\ds \limt |\vla(t) |_0 = 0 ,
}
then we can obtain a contradiction.

For any $\de > 0$, define $g_\de \in C^0(\bOm)$ by
$$
g_\de := g(x,\de) , \quad x \in \bOm .
$$
By the properties of $g$, and the principal eigenvalue function $\mu_0(\cdot)$ (see Lemma~\ref{princ_eval.lem} and \cite{CUE}),
we have
$$
\text{ $ g_\de \le g_0$ \ and \ $\lim_{\de \todowns 0} | g_\de - g_0 |_0 = 0$}
\implies
\text{ $ \mu_0(g_\de) \ge \mu_0(g_0) $ \ and \ $\lim_{\de \todowns 0} \mu_0(g_\de) = \mu_0(g_0) $}
$$
(the final limiting result is not explicitly stated in \cite{CUE}, but it can readily be proved using the minimisation characterisation of $\mu_0(\rho)$ in $(1.3)$ of \cite{CUE}; the argument is similar to the proof of \cite[Proposition~4.3]{CUE}).
Hence, since $ \la > \lamin(g) = \mu_0(g_0) $,
we may choose $\de $ sufficiently small that $ \la > \mu_0(g_\de) $.
Also, by \eqref{vla_to_z.eq}, we may choose $t_\de \ge 0$ such that
\eqn{  \label{vla_small_tde.eq}
| \vla(t) |_0 \le \de/4 , \quad t \ge t_\de  .
}
Now, by regarding $t_\de$ as the initial time, and $ v_\de := \vla(t_\de)$ as the initial value, we can follow the argument in Section~\ref{proof-part_c.sec} to show that
$ T_{\la g, v_\de} \le T_{\la g_\de, v_\de} < \infty , $
that is, $ v_{\la g, v_\de} $ blows up in finite time
(the inequality $ \la > \mu_0(g_\de) $ provides the analogue here of the
inequality $\la > \lamax(g) = \mu_0(g_\infty)$ used in Section~\ref{proof-part_c.sec}).
This clearly contradicts \eqref{vla_small_tde.eq}, and so proves Lemma~\ref{vinfty_ne_0.lem}.
\eop

\end{document}